\pdfoutput=1
\documentclass[12pt]{notices}
\usepackage[dvipsnames]{xcolor}
\usepackage{amsmath, amsthm, amssymb, hyperref, color}
\usepackage{mathtools}
\hypersetup{
     colorlinks = true,
     linkcolor = black,
     anchorcolor = black,
     citecolor = black,
     filecolor = black,
     urlcolor = black
     }
\usepackage{graphicx}
\usepackage[font=small]{caption}
\usepackage{cleveref}
\usepackage{enumerate}
\usepackage{caption}
\usepackage[all]{xypic}
\usepackage{verbatim}
\usepackage{chemfig, chemnum}
\usepackage{tikz}
\usetikzlibrary{decorations.markings}
\usepackage[linesnumbered,lined,commentsnumbered]{algorithm2e}
\usepackage{caption}
\usepackage{subcaption}
\oddsidemargin \evensidemargin
\usepackage{makecell}
\usepackage{array}
\newcolumntype{?}{!{\vrule width 1pt}}

\newtheorem{theorem}{Theorem}

\newtheorem{example}[theorem]{Example}

\theoremstyle{definition}

\newcommand{\PP}{\mathbb{P}}
\newcommand{\RR}{\mathbb{R}}

\newcommand{\CC}{\mathbb{C}}
\newcommand{\ZZ}{\mathbb{Z}}
\newcommand{\NN}{\mathbb{N}}

\renewcommand{\d}{\mathrm{d}}

\newcommand{\cI}{\mathcal{I}}

\DeclareMathOperator{\Gr}{Gr}
\DeclareMathOperator{\dlog}{dlog}

\usetikzlibrary{snakes}
\usetikzlibrary{calc}
\usetikzlibrary{decorations.pathmorphing}

\tikzset{
    quark/.style={
        decoration={},
        decorate
    },
    lepton/.style={
        decoration={},
        decorate
    },
    gluon/.style={
        decoration={coil, aspect=0.75, mirror, segment length=1.5mm},
        decorate
    },
    gluon_small/.style={
        decoration={coil, aspect=0.75, mirror, segment length=1mm, amplitude =0.6mm},
        decorate
    },
    vector/.style={
        decoration={snake, aspect=0.75, mirror, segment length=2mm},
        decorate
    },
    Higgs/.style={
        decoration={},
        densely dashed,
        decorate
    },
}

\title{\bf Algebraic \& Positive Geometry of the Universe:\\from Particles to Galaxies}

\author{\bf Claudia Fevola and Anna-Laura Sattelberger}
\date{}

\begin{document}
\maketitle
\begin{abstract}
In recent years, the intersection of algebra, geometry, and combinatorics with particle physics and cosmology has led to significant advances. Central to this progress is the twofold formulation of the study of particle interactions and observables in the universe: on the one hand, Feynman's approach reduces to the study of intricate integrals; on the other hand, one encounters the study of positive geometries. This article introduces key developments, mathematical tools, and the connections that drive progress at the frontier between algebraic geometry, the theory of $D$-modules, combinatorics, and physics. All these threads contribute to shaping the flourishing field of positive geometry, which aims to establish a unifying mathematical language for describing phenomena in cosmology and particle physics.
\end{abstract}

\section{Introduction}\label{sec:intro}
The relationship between mathematics and physics is profoundly symbiotic: mathematics provides the language and tools to describe and predict physical phenomena, while physics inspires the discovery and development of new mathematical concepts. 
Many significant advances in physics rely on 
sophisticated mathematical frameworks. Historical examples of this interplay abound. This dynamic relationship continues today in fields such as quantum field theory and cosmology, where novel mathematical concepts and fundamental physics are intricately intertwined. 

In this article, we explore the recent developments at the intersection of algebraic geometry, algebraic analysis, combinatorics, and their applications to particle physics and cosmology. We highlight key research directions, innovative ideas, and the mathematical structures that have emerged or have been revitalized in response to challenges arising in physics.

\begin{figure}[h]
\centering
\tikzset{every picture/.style={line width=0.75pt}} 

\begin{tikzpicture}[x=0.75pt,y=0.75pt,yscale=-1,xscale=1]

\draw  [fill={rgb, 255:red, 189; green, 185; blue, 187 }  ,fill opacity=.7][line width=1.3]  (101,135) .. controls (101,121.19) and (112.19,110) .. (126,110) .. controls (139.81,110) and (151,121.19) .. (151,135) .. controls (151,148.81) and (139.81,160) .. (126,160) .. controls (112.19,160) and (101,148.81) .. (101,135) -- cycle ;
\draw [line width=1.3]    (69.5,99.64) -- (100.23,121.33) ;
\draw [shift={(103.5,123.64)}, rotate = 215.22] [fill={rgb, 255:red, 0; green, 0; blue, 0 }  ][line width=0.08]  [draw opacity=0] (11.61,-5.58) -- (0,0) -- (11.61,5.58) -- cycle    ;
\draw [line width=1.3]    (70.5,170.64) -- (100.27,148.99) ;
\draw [shift={(103.5,146.64)}, rotate = 143.97] [fill={rgb, 255:red, 0; green, 0; blue, 0 }  ][line width=0.08]  [draw opacity=0] (11.61,-5.58) -- (0,0) -- (11.61,5.58) -- cycle    ;
\draw [line width=1.3]    (149.5,126.64) -- (182.71,115.42) ;
\draw [shift={(186.5,114.14)}, rotate = 161.33] [fill={rgb, 255:red, 0; green, 0; blue, 0 }  ][line width=0.08]  [draw opacity=0] (11.61,-5.58) -- (0,0) -- (11.61,5.58) -- cycle    ;
\draw [line width=1.3]    (150,140) -- (185.5,141.03) ;
\draw [shift={(189.5,141.14)}, rotate = 181.65] [fill={rgb, 255:red, 0; green, 0; blue, 0 }  ][line width=0.08]  [draw opacity=0] (11.61,-5.58) -- (0,0) -- (11.61,5.58) -- cycle    ;
\draw [line width=1.3]    (61.5,134.64) -- (97,134.96) ;
\draw [shift={(101,135)}, rotate = 180.52] [fill={rgb, 255:red, 0; green, 0; blue, 0 }  ][line width=0.08]  [draw opacity=0] (11.61,-5.58) -- (0,0) -- (11.61,5.58) -- cycle    ;
\draw [line width=1.3]    (145.5,151.64) -- (181.74,164.78) ;
\draw [shift={(185.5,166.14)}, rotate = 199.93] [fill={rgb, 255:red, 0; green, 0; blue, 0 }  ][line width=0.08]  [draw opacity=0] (11.61,-5.58) -- (0,0) -- (11.61,5.58) -- cycle    ;

\draw (117,127.4) node [anchor=north west][inner sep=0.75pt]  [font=\Large]  {$\mathcal{S}$};
\draw (102,160.54) node [anchor=north west][inner sep=0.75pt]  [font=\LARGE]  {$.$};
\draw (128,165.54) node [anchor=north west][inner sep=0.75pt]  [font=\LARGE]  {$.$};
\draw (140,160.54) node [anchor=north west][inner sep=0.75pt]  [font=\LARGE]  {$.$};
\draw (114,166.54) node [anchor=north west][inner sep=0.75pt]  [font=\LARGE]  {$.$};
\draw (52,102.54) node [anchor=north west][inner sep=0.75pt]    {$p_{1}$};
\draw (48,137.54) node [anchor=north west][inner sep=0.75pt]    {$p_{2}$};
\draw (52,168.54) node [anchor=north west][inner sep=0.75pt]    {$p_{3}$};
\draw (190,137.54) node [anchor=north west][inner sep=0.75pt]    {$p_{n-1}$};
\draw (191,102.54) node [anchor=north west][inner sep=0.75pt]    {$p_{n}$};
\draw (187.5,168.54) node [anchor=north west][inner sep=0.75pt]    {$p_{n-2}$};
\end{tikzpicture}
\caption{A graphical representation of a scattering process between $n$ particles, such 
as happens for instance in a particle accelerator.
}
\label{fig:scattering_process}
\end{figure}
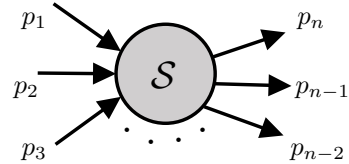

Particle physics studies the interactions of elementary particles that constitute matter and radiation. Mathematically, scattering amplitudes measure the probabilities of specific outcomes in these interactions,
which are studied in practice in particle accelerators by experimentalists.
\Cref{fig:scattering_process} provides a pictorial representation of a scattering process involving $n$ particles. The gray disc represents the scattering of the particles in the process. The goal is to predict what happened during this interaction by analyzing the behavior of the outgoing particles. By doing so, physicists hope to answer some of the most fundamental questions about the nature of the universe, such as how it began, what it is made of, and how it will evolve in the future.

More formally, the scattering amplitude associated to a scattering process is a complex-valued function
\begin{align}\label{eq:scattering_amplitude}
\left( \RR^{1,D-1} \right)^n \to \CC ,\,\left(p_1,\ldots,p_n\right)\mapsto A\left(p_1,\ldots,p_n\right),
\end{align}
where the $p_i = (p_i^0, p_i^1,\dots, p_i^{D-1})$ for $i = 1,2,\dots,n$ are the {momentum vectors} associated to the particles. These are elements of the $D$-dimensional {\em Minkowski space} $\RR^{1,D-1}$, i.e., $\RR^D$ endowed with the indefinite inner product $p\cdot q = p^0q^0 - p^1q^1 -\cdots - p^{D-1}q^{D-1}$. We write $p^2$ to denote $p\cdot p$. The most relevant case in physics is when $D=4$, i.e., three space-dimensions and one dimension representing time. Momentum vectors capture physical information about the particles in the scattering experiment. The modulus squared of the scattering amplitude \eqref{eq:scattering_amplitude} can be interpreted as a joint probability density function describing what to expect for the outcome of the experiment. 

In perturbation theory, there exist two primary ways to compute scattering amplitudes:
\begin{enumerate}
\item \textbf{Feynman integrals}
A breakthrough in the computation of amplitudes came in the 1940's with Richard Feynman’s diagrammatic method. Feynman diagrams encode the dynamics of particle collisions, with explicit rules for associating a function to each diagram, which calculates the process’s probability. 
The scattering amplitude is expressed as a sum over all possible Feynman diagrams~$G$:  
\[A \,=\, \sum_{G} a_G \cdot \cI_G \, .\] 
Here, $\cI_G$ represents the Feynman integral~\cite{Weinzierl} associated with a graph $G$ (and a rational function when $G$ is a tree), and it is defined according to Feynman rules. 
The coefficient $a_G$ depends both on the properties of the particles involved and on the topology of the graph. Studying the analytic properties of the amplitude $A$ often involves analyzing the terms~$\cI_G$ individually. The analytic and numerical evaluation of Feynman integrals, as well as the investigation of their monodromy representations and branch cuts as multi-valued functions, are difficult problems that have sparked numerous lines of research. For example, identifying the special functions required to express these integrals has been a difficult question studied since the early days of quantum field theory (QFT) and remains an active field of research. 

\item \textbf{Positive geometries\footnote{Not to be confused with the field of positive geometry~\cite{LeMatematiche},  
where positive geometries are one key topic among others.}} 
A more recent perspective avoids the splitting into individual Feynman diagrams. It defines scattering amplitudes in terms of purely geometric structures---depending only on the initial and final particle states, without referencing their evolution in spacetime. The need for a reformulation of the scattering amplitude inspired Arkani-Hamed and Trnka to introduce the amplituhedron~\cite{Amplitu} for $\mathcal{N}=4$ supersymmetric Yang--Mills (SYM) in the planar limit, a certain QFT. 
The fundamental idea behind its definition is to provide a mathematical object whose volume directly computes the scattering amplitude. This motivation has later inspired the more general and rigid notion of positive geometries, which we address in~\Cref{sec:combalggeo}. In the conception of positive geometry, the connection between the geometry and physics is provided by a differential form, called ``canonical form.'' This form is associated to a semi-algebraic subset of the real points of a variety, or to a subvariety in more recent literature~\cite{BrownDupont}. It is uniquely determined by the requirement of having at most simple poles on all boundary components of the space, with residues along each boundary component given by the canonical form of that component. The mathematical properties of amplituhedra and positive geometries exhibit rich structures, suggesting that they should be studied \hbox{as mathematical concepts on their own.}
\end{enumerate}

In theoretical cosmology, Feynman's approach and positive geometries are  highly relevant as well.
A fundamental aspect of this branch of physics is to understand how matter and energy are distributed throughout the universe at any given time. In particle physics, scattering experiments are run in particle accelerators---such as in the world's largest and most powerful one, the Large Hadron Collider 
at CERN. However, one could argue that the Big Bang already did such experiments for us and that the resulting data are out there in the~sky. To extract them, new mathematics is needed. One approach for addressing this challenge is to analyze large-scale correlations observed today, such as among galaxies, and use them to infer the physical processes that shaped the early universe. Simultaneously, as a theoretical foundation, cosmologists rely on frameworks grounded in quantum field theory and \mbox{positive geometry}.

\begin{figure}[h]
\begin{center}
\includegraphics[width=5.8cm]{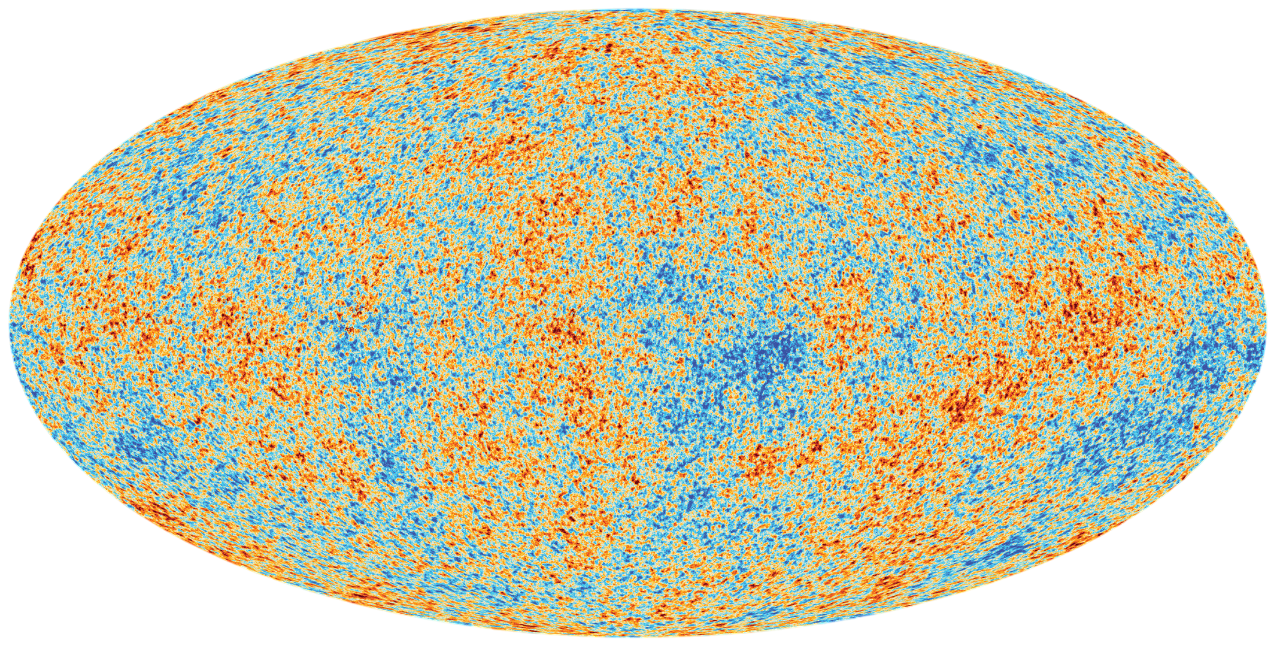}
\end{center}
\caption{The CMB as observed by Planck. Credit: \href{https://www.esa.int/ESA_Multimedia/Images/2013/10/The_Cosmic_microwave_background_CMB_as_observed_by_Planck}{ESA, Planck Collaboration}, 2013.}
\label{fig:CMB}
\end{figure}

The cosmic microwave background (CMB), shown in \Cref{fig:CMB}, is a remnant of the first light in our universe and is evidence that the Big Bang happened approximately $13{.}8$ billion years ago.

\vspace*{2mm}

The statistical properties of the initial conditions of our universe are modeled by ``cosmological correlators.'' Despite foundational differences in the principles used to model scattering amplitudes and theoretical cosmology, interestingly, the methods for computing these correlators rely on techniques developed for studying Feynman integrals. In fact, in the perturbative approach, correlators are also associated with Feynman diagrams. In the context of cosmological correlators, already trees---i.e., graphs without cycles---are significantly intricate to study; however, recent advances~\cite{DEcosmological} revealed that the resulting combinatorial structures are unexpectedly well-behaved. The mathematics of cosmological correlators is a relatively recent area of study. In particular, the interplay with combinatorics is novel. There are many intriguing and challenging mathematical problems that are awaiting to~be~tackled.

Mathematically, both Feynman integrals and correlation functions are inherently linked to generalized Euler integrals~\cite{GKZ90}, i.e., integrals of the form
\begin{align}\label{eq:integral}
\mathcal{I}(c) =\! \int_\Gamma f_1^{s_1}\cdots f_k^{s_k}x_1^{\nu_1}\cdots x_n^{\nu_n} \, \frac{\mathrm{d}x_1}{x_1}\wedge \cdots \wedge \frac{\mathrm{d}x_n}{x_n}, 
\end{align}
where $f_1,\ldots,f_k\in \CC[x_1^{\pm 1},\ldots,x_n^{\pm 1}]$ are Laurent polynomials, the exponents $s_i$ and $\nu_j$ are complex numbers, and $\Gamma$ is a suitable integration contour. The integral~\eqref{eq:integral} is to be read as a function of the coefficients of the Laurent polynomials $f_j=\sum_{u \in \ZZ^n} c^{(j)}_u x^u$. In particle physics, $k=1$, and $f$ is the graph polynomial of a Feynman diagram. In cosmology, for certain toy models, 
the $f_j$'s define hyperplanes in~$\CC^n$, the $s_j$'s are integers, 
and the $\nu_i$'s (later on denoted by~$\epsilon_i)$ are related to the assumed model of the universe. Diverse mathematical tools are needed to study these integrals---reaching from analysis, algebraic geometry, mixed Hodge theory, analytic number theory, intersection theory, to numerical evaluation methods. This list still does not cover the full range of areas involved; we here focus on a selection that is 
biased towards~algebra. 

In this introductory exposition, we present the novel and mutually enriching interactions arising from employing algebra, geometry, and combinatorics as tools for elaborating the mathematical concepts designed by theoretical physicists to investigate the universe. Algebraic geometry as a unifying language allows us treating physical phenomena on both small scales (interactions of fundamental particles) and large scales (structure of the universe).~Vice versa, physics inspires new mathematics, such as the flourishing field of positive geometry, which aims to establish a new mathematical field as a common language and to formulate the concept of positivity for some mathematical objects. 

\vspace{-2.5mm}
\paragraph{Timeliness}
The topicality and the active movement of the subject is, for instance, mirrored by the two recently awarded ERC synergy grants ``\href{https://positive-geometry.com/}{UNIVERSE+}: Positive Geometry in Cosmology and Particle Physics'' led by Nima Arkani-Hamed, Daniel Baumann, Johannes Henn, and Bernd Sturmfels, and ``MaScAmp: Mathematics of Scattering Amplitudes'' by Ruth Britto, Francis Brown, Axel Kleinschmidt, and Oliver Schlotterer. The research community is looking forward to a variety of upcoming interdisciplinary events within the next years. These will be rooted throughout Europe as well as in the US, with the IAS Princeton being one of the main nodes.

\vspace{-3mm}
\paragraph{Outline} We divide our presentation into four sections: we address particle physics, cosmology, the mathematical toolkit (algebraic geometry, algebraic analysis, and combinatorial algebraic geometry), and current interactions among these fields. 

\section{Particle physics}\label{sec:partphy}
\vspace*{-1mm}%
This section provides a self-contained introduction to Feynman integrals in particle physics~\cite{Weinzierl}, specifically their close bond with generalized Euler integrals~\cite{GKZ90}.

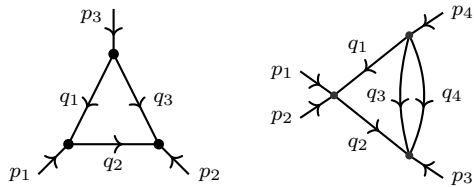
\begin{figure}[h]
    \centering
    \begin{subfigure}[c]{0.2\textwidth}
     \centering
\begin{tikzpicture}[line width=0.7,scale=1.0]
\begin{scope}[very thick,decoration={
    markings,
    mark=at position 0.6 with {\arrow{>}}}]
    \coordinate (v1) at (0,0);
    \coordinate (v2) at (0.6,1.2);
    \coordinate (v4) at (1.2,0);

    \coordinate (v11) at (-0.4,-0.4);
    \coordinate (v21) at (1.6,-0.4);
    \coordinate (v41) at (0.6,1.8);
    \draw[thick,postaction={decorate}] (v2) -- (v1);
    \draw[thick,postaction={decorate}] (v2) -- (v4);
    \draw[thick,postaction={decorate}] (v1) -- (v4);
    \draw[thick,postaction={decorate}] (v11) -- (v1);
    \draw[thick,postaction={decorate}] (v41) -- (v2);
    \draw[thick,postaction={decorate}] (v21) -- (v4);
    \foreach \point in {v1, v2, v4} {
        \fill[black] (\point) circle [radius=0.07];
    }

	\node at (1.55,0.6) [font=\footnotesize, anchor=east] {$q_3$};
	\node at (0.3,0.6) [font=\footnotesize, anchor=east] {$q_1$};
	\node at (0.6,-0.05) [font=\footnotesize, anchor=north] {$q_2$};
	\node at (-0.35,-0.4)
    [font=\footnotesize, anchor=east] {$p_1$};
	\node at (0.6,1.7)
    [font=\footnotesize, anchor=east] {$p_3$};
	\node at (v21) 
    [font=\footnotesize, anchor=west] {$p_2$};
    \end{scope}
    \end{tikzpicture}
    \end{subfigure}
    \begin{subfigure}[c]{0.2\textwidth}
\centering
\begin{tikzpicture}[line width=1,scale=1.0]
\begin{scope}[very thick,decoration={
    markings,
    mark=at position 0.6 with {\arrow{>}}}]
    \coordinate (v1) at (-1,0);
    \coordinate (v2) at (0,0.8);
    \coordinate (v3) at (0,-0.8);
    \draw[thick,postaction={decorate}] (v2) -- (v1);
    \draw[thick,postaction={decorate}] (v1) -- (v3);
    \draw[thick,postaction={decorate}] (v2) to[out=-65,in=65] (v3);
    \draw[thick,postaction={decorate}] (v2) to[out=-105,in=105] (v3);
    \draw[thick,postaction={decorate}] (-1.45,-0.3) -- (v1);
    \draw[thick,postaction={decorate}] (-1.45,0.32) -- (v1);
    \draw[thick,postaction={decorate}] (0.45,1.1) -- (v2);
    \draw[thick,postaction={decorate}] (0.45,-1.1) -- (v3);
    \foreach \point in {v1, v2, v3} {
    	\fill[white!20!black] (\point) circle (1.5pt);
    }
	\node at (-0.4,-0.65) [font=\footnotesize, anchor=east] {$q_2$};
	\node at (-0.15,0) [font=\footnotesize, anchor=east] {$q_3$};
	\node at (-0.4,0.65) [font=\footnotesize, anchor=east] {$q_1$};
	\node at (0.85,0) [font=\footnotesize, anchor=east] {$q_4$};
	\node at (-1.4,0.3) [font=\footnotesize, anchor=east] {$p_1$};
	\node at (-1.4,-0.3) [font=\footnotesize, anchor=east] {$p_2$};
	\node at (1,1.1) [font=\footnotesize, anchor=east] {$p_4$};
	\node at (1,-1.1) [font=\footnotesize, anchor=east] {$p_3$};
    \end{scope}
\end{tikzpicture}
\end{subfigure}  
\caption{Some Feynman diagrams: the ``one-loop triangle'' (left) and the ``parachute'' (right) diagram.
}
\label{fig:Feynman}
\end{figure}

A Feynman diagram $G$ is a connected, oriented graph $(V,E)$ with $n$ ``external legs'' (open edges) attached, $m$ internal edges, $|V|$ vertices, and $\ell$ independent cycles, where $\ell = m - |V| + 1,$ since $G$ is connected. The external legs represent the incoming and outgoing particles, labeled with a momentum vector $p_i \in \RR^{1,D-1}$ each. 
Each edge $e_i$ in $G$ is arbitrarily oriented, but the associated integral will be independent of this choice. Each internal edge instead is equipped with a momentum vector $q_i \in \RR^{1,D-1}$ and a mass $m_i \in \RR_{\geq 0}$ of the intermediate particle propagating along the edge~$e_i$. Some examples of Feynman diagrams are depicted in \Cref{fig:Feynman}. 
Momentum conservation at each vertex states that the sum of incoming momenta equals the sum of outgoing ones, yielding a linear equation in $p_i$ and $q_i$ for each~vertex.
To each such diagram, one associates a Feynman integral~$\cI_G$.

There are several representations of Feynman integrals in common use, which can be transformed into each other by changes of variables. For example, some representations use variables associated to cycles, namely the loop momenta, as integration variables, while others replace these with auxiliary \mbox{variables} known as ``Schwinger parameters.'' 

Here, we focus on the Lee--Pomeransky representation, as it is most directly related to generalized Euler integrals~\eqref{eq:integral}. To each internal edge $e_i$, we associate a Schwinger parameter \( \alpha_i \in {\mathbb{C}}^* \). The graph polynomial will be a polynomial in the $\alpha_i$'s. To introduce the Lee--Pomeransky representation of a Feynman integral, we first introduce some necessary graph-theoretical concepts. A {\em subtree} of a graph~$G$ is a connected, acyclic subgraph of~$G$. A {\em spanning tree} of $G$ is a subtree that contains all vertices of~$G$.
We write $\mathcal{T}_1$ for the set of all spanning trees of~$G$. 

The {\em first Symanzik polynomial} of $G$ is 
$$\mathcal{U}_G\, \coloneqq\, \sum_{T\in \mathcal{T}_1}\prod_{e_i\notin T} \alpha_i\, ,$$
where the product runs over the $L$ internal edges that are removed from $G$ to obtain the spanning tree $T$.

A {\em spanning $2$-forest} $T'$ of $G$ is a disjoint union of two subtrees of $G$ whose union contains all vertices of $G$. We write $T_1$ and $T_2$ for its connected components, and denote by $\mathcal{P}_{T_i}$ the set of external momenta attached to~$T_i$. 

The {\em second Symanzik polynomial}, denoted $\mathcal{F}_G$, is
\begin{equation*}
\sum_{T'\in \mathcal{T}_2} \! \Bigl(\sum_{p_j\in \mathcal{P}_{T_1}}\sum_{p_k\in \mathcal{P}_{T_2}} \frac{p_j\cdot p_k}{\mu^2}\Bigr)  \! \prod_{e_i\not\in T'} \alpha_i - \,  \mathcal{U}_G\cdot \! \sum_{e_i\in E} \frac{m_i^2}{\mu^2} \alpha_i \, ,
\end{equation*}
where $\mu$ is a scaling factor to make the polynomial dimensionless, that we will ignore henceforth. The \hbox{$p_j\cdot p_k$} and $m_i^2$ are called ``kinematic variables.''

Note that the first and second Symanzik polynomial are homogeneous polynomials in the Schwinger parameters of degree $\ell$ and $\ell + 1$, and at most linear and quadratic in each $\alpha_i$, respectively. Their sum, $\mathcal{G}_G=\mathcal{U}_G+\mathcal{F}_G$, is the \textit{graph polynomial} of~$G$.

\begin{example}\label{ex:FeynmanInt}
For $G$ the parachute Feynman diagram from \Cref{fig:Feynman}, the Symanzik polynomials are
\begin{align*}
\mathcal{U}_G  & \,=\, 
\alpha_1 \alpha_3 + \alpha_2 \alpha_3 + \alpha_1 \alpha_4 + \alpha_2 \alpha_4 + \alpha_3 \alpha_4\, , \\
\mathcal{F}_G & \,=\,\left(p_1+p_2\right)^2\alpha_1\alpha_2 \left(\alpha_3+\alpha_4\right) + p_3^2\alpha_2\alpha_3\alpha_4 \\ &\quad \ +p_4^2\alpha_1\alpha_3\alpha_4 - \mathcal{U}_G \cdot \big(\textstyle{\sum}_{i=1}^4m_i^2\alpha_i\big) \, .
\end{align*}
The algebraic variety of $\mathcal{G}_G$ is depicted in \Cref{fig:var}.
\end{example}

Whenever $\ell>1$, the {\em Feynman integral} of $G$ (in the Lee--Pomeransky representation) is
\begin{align*}\label{eq:Lee_Pomeranski}
\cI_G \,=\, N_\nu\cdot \int_{\RR_{>0}^n} \big( \prod_{i = 1}^n \alpha_i^{\nu_i-1}\big)\mathcal{G}_G^{-D/2} \,  \d \alpha \, ,
\end{align*}
where $N_\nu$ is a normalizing factor that depends on the integer exponents $\nu = (\nu_1,\dots,\nu_n)$, and $\d \alpha$ denotes the $n$-form $\d \alpha_1\wedge\cdots\wedge \d \alpha_n$.

\section{Cosmology}\label{sec:cosmo}
Feynman diagrams also come into play when studying the initial conditions of our universe. In this context, graphs are considered as undirected and, in order to study cosmological integrals, one truncates the external edges for simplicity. Therefore, in this section, a Feynman diagram $G$ is a connected, undirected graph without external legs. We denote by $V=V_G$ the set of its $n$ vertices $v_i$, and  $E=E_G$ is the collection of its edges $ij$ for some $v_i,v_j\in V$. The vector 
$(X,Y) \coloneqq (X_1,\dots,X_n,\{Y_{ij}\}_{ij\in E})$ represents the energies associated to the vertices and edges of the graph. The $X_i$'s represent external outgoing energies, while the $Y_{ij}$'s are energies exchanged in the process and therefore are real and positive numbers, although it is useful to think of them in larger domains, using analytic continuation.

Cosmological integrals are defined by means of the {flat space wavefunction coefficient} $\psi_{\text{flat}}^G$ (or simply~$\psi_{\text{flat}}$, when the graph is clear from the context). Concretely, $\psi_{\text{flat}}^G$ is a rational function of the parameters $(X,Y)$ that can be computed in three ways:
\begin{enumerate}[(1)]
\item using the Feynman rules, as in \cite[Section~2.2]{DEcosmological},
\item via the recursion formula 
from~\cite[Section~2]{cosmowave}, or
\item as the rational function defining the canonical form of a cosmological polytope, see~\cite[(3.5)]{cosmowave}.
\end{enumerate}
It is a conjecture in physics that for any Feynman diagram, all these methods lead to the same rational function. Here, we focus on the third approach.

The cosmological polytope $P_G$ of a Feynman diagram $G$ plays the role of a positive geometry, which we will encounter in detail in \Cref{sec:combalggeo}. Given a Feynman diagram $G=(V,E)$, the cosmological polytope $P_G$ lives in the space $\RR^{|V|+|E|}$ with basis $(X,Y) \coloneqq (X_1,\dots,X_n,\{Y_{ij}\}_{ij\in E})$. The facets of $P_G$ are in bijection with the connected subgraphs of $G$. In practice, given a connected subgraph $H = (V_H,E_H)$ of~$G$, the facet $F_H$ is the intersection of $P_G$ with the hyperplane $L_H$ defined by the polynomial
\begin{equation*}\label{eq:linear_forms_subgraphs}
\sum_{v_i \,\in \, V_H} \! X_i \ \ \  + \! \sum_{\substack{e \,=\, ij \in V,\\ \,v_i\,\in\, V_H, \,v_j \,\not\in\, V_H}} \!  Y_{ij} \ \ +\! 
    \sum_{\substack{e \,=\, ij \,\not\in\, E_H,\\ v_i,v_j\,\in\, V_H}} 2Y_{ij} \, ,
\end{equation*}
see~\cite{cosmowave} for references.
Polytopes are indeed examples of positive geometries~\cite{PosGeom, BrownDupont}. The condition on the singularities required in the definition of the canonical form guarantees that it can be written as
\begin{equation}\label{eq:canonical_form_cosmo}
\frac{P}{\prod_{H\subseteq G} L_H} \,\d X\wedge\d Y  ,
\end{equation}    
where $\d X = \d X_1 \wedge \cdots \wedge \d
X_n$, $\d Y = \wedge_{ij\in E}\, \d Y_{ij}$, and the product in the denominator runs over all connected subgraphs of~$G$.
The numerator $P\in \mathbb{C}[X,Y]$ defines the adjoint hypersurface, which geometrically is the hypersurface of minimal degree that contains all linear spaces that are intersections of hyperplanes $L_H$ and do not contain any face of $P_G$. The intuition behind the definition of the cosmological polytope $P_G$ is to build a polytope whose canonical form~\eqref{eq:canonical_form_cosmo} equals 
\[\psi_{\text{flat}}^G(X,Y) \, \d X \wedge\d Y  .\]  
In a certain cosmological toy model, see~\cite{DEcosmological,cosmowave}, the {\em cosmological integral} now is defined to~be
\begin{equation}\label{eq:cosmoint_explicit}
    \int_{\RR^n_{>0}} 2^{n-1}\bigl(\prod_{ij\in E} Y_{ij}\bigr)\cdot\widetilde{\psi}_{\text{flat}}\cdot \alpha_1^{\epsilon_1}\cdots \alpha_n^{\epsilon_n} \, \d\alpha \,,
\end{equation}  
where $\d \alpha  \coloneqq \d \alpha_1\wedge \dots \wedge \d\alpha_n,$ and $\widetilde{\psi}_{\text{flat}}$ is obtained from $\psi_{\text{flat}}$ by shifting the 
$X_i$'s by the variables~$\alpha_i$, i.e., replacing $X_i\to X_i+\alpha_i$. Physically, the role of the integration variables $\alpha$ is to parameterize the time in which a physical process happens, while the
exponents $\epsilon_i$ are parameters of the cosmology in which the process takes place. In simple terms, a cosmology refers to the overall structure and behavior of the universe that serves as the backdrop for the physical process being discussed. 
Some integer values of $\epsilon_i$  correspond to significant cases of cosmological models of the universe, see~\cite[p.~20]{DEcosmological}.
However, also non-integer values of $\epsilon_i$ are interesting for physics, so we consider the $\epsilon_i$ to be generic, complex parameters.
\begin{example}\label{ex:2-site}
Figure~\ref{fig:2-site} depicts  the path graph on two vertices, referred to as the ``$\,2$-site chain.'' Its associated cosmological polytope lives in the hyperplane $\{X_1+X_2+Y_{12}=1\}$ in $\RR^3$ and its facet hyperplanes $\{X_1+X_2=0\}$, $\{X_1+Y_{12}=0\}$, $\{X_2+Y_{12}=0\}$ can be obtained via the formula in~\eqref{eq:linear_forms_subgraphs}. As these hyperplanes intersect only in the vertices of the polytope, their adjoint is $P=1$ so that 
$$\psi_{\text{flat}} \, = \, \frac{1}{(X_1+X_2)(X_1+Y_{12})(X_2+Y_{12})} \, .$$

\begin{figure}
\centering
    \begin{tikzpicture}
     \draw (-4.5,0) -- (-3,0);
     \filldraw[Aquamarine] (-4.5,0) circle (3pt); 
     \filldraw[Bittersweet] (-3,0) circle (3pt); 
     \node at (-4.5,-.4) {\small$X_1$};
    \node at (-3,-.4) {\small$X_2$};
    \node at (-3.75,.25) {\small$Y_{12}$};
     \draw (-1,-.7) -- (1,-.7) -- (0,.7) -- (-1,-.7);
    \filldraw[Aquamarine] (-.8,0) circle (1.8pt);
    \filldraw[Bittersweet] (.8,0) circle (1.8pt);
    \filldraw[black] (-.3,-.94) circle (1.7pt);
    \filldraw[black] (.3,-.94) circle (1.7pt);
    \draw (-.3,-.94) -- (.3,-.94);
    \node at (-1.5,-1) {\small $(-1,1,1)$};
    \node at (1.5,-1) {\small $(1,-1,1)$};
    \node at (0,.95) {\small $(1,1,-1)$};
    \end{tikzpicture}
    \caption{The 2-site chain graph (left) and its associated cosmological polytope (right).}
\label{fig:2-site}
\end{figure}
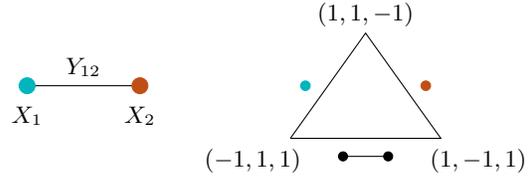

\noindent The associated cosmological integral is derived from~\eqref{eq:cosmoint_explicit}, where 
$\widetilde{\psi}_{\text{flat}}$ equals
\begin{small}
\[\frac{1}{(\alpha_1+\alpha_2+X_1+X_2)(\alpha_1+X_1+Y_{12})(\alpha_2+X_2+Y_{12})} \, .\]
\end{small}
\end{example}

One may ask why these integrals are relevant for mathematics. They serve as meaningful examples of generalized Euler integrals \eqref{eq:integral}, where all the polynomials $f_i$ have degree~$1$.
Many properties of such families of integrals can be derived from studying the complement in $(\CC^*)^n$ of the hyperplanes defined by the polynomials in the integral. This locus will be introduced for polynomials of higher degree in \mbox{\Cref{sec:alggeo}}. 

The case of hyperplanes that are in general position has mathematically been extensively studied by Aomoto, Kita, Orlik, and Terao. 
However, in the context of cosmological integrals, the coefficients of the polynomials are not arbitrary. Instead, they are constrained to lie in specific linear subspaces determined by the physical energies of the system. This restriction makes these integrals particularly special. For fixed values of $X_1,X_2,Y_{12}$, the complement of hyperplanes from \Cref{ex:2-site} lives in the torus $(\CC^*)^2$ with coordinates $(\alpha_1,\alpha_2)$, and it is represented in~\cite[Figure~1]{DEcosmological}. Investigating the singularities of cosmological integrals and the differential equations that they satisfy---ultimately aiming to determine closed-form expressions of the resulting functions---requires the utilization of a wide range of mathematical tools.

\section{Mathematical toolkit}\label{sec:mathtools}
Theoretical physics and cosmology increasingly benefit from methods from the fields of algebra and combinatorics,
whose pure and rigorous structures naturally align with the analysis of fundamental interactions. For instance, the integrals defined in~\eqref{eq:integral} are multi-valued functions. They have singularities, the allowance of non-generic parameters is needed, and a high-precision evaluation of them is required to test physical theories. To tackle these challenges, we present tools from algebraic geometry and nonlinear algebra~\cite{NonlinearAlgebra}, the theory of \mbox{$D$-modules}~\cite{HTT08}, as well as combinatorial approaches~\cite{Lam,Williams}.

\subsection{Algebraic geometry}\label{sec:alggeo}
\begin{figure}[h]
\centering
\includegraphics[width=5.2cm]{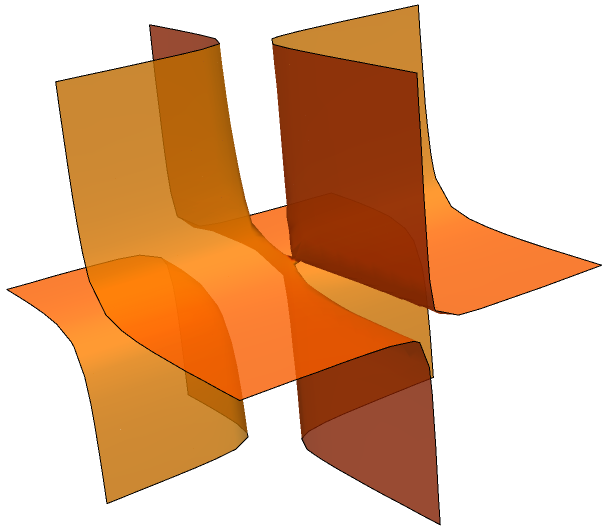}
\caption{
The set of real zeros of the graph polynomial~$\mathcal{G}_G$, where $G$ is the massless parachute diagram as in \Cref{ex:FeynmanInt}, with parameters $(p_1+p_2)^2=25,$ $p_3^2=49,$ $p_4^2=9$.
}
\label{fig:var}
\end{figure}

Affine algebraic varieties $V(I)\subset \CC^n$, with $I\subset \CC[x_1,\ldots,x_n]$, encode common vanishing sets of systems of polynomials. The underlying topology is the Zariski topology, which is defined in terms of prime ideals and is coarser than the Euclidean topology. An example of a variety is depicted in \Cref{fig:var}.

\vspace*{1mm}

We denote the affine $n$-space 
by its closed points, i.e., $\CC^n$. 
In the same vein, we denote the algebraic $n$-torus 
short-hand by~$(\CC^\ast)^n$. Varieties that allow for a closed embedding into an algebraic torus are called ``very affine.'' Examples include complements of hyperplane arrangements in an algebraic torus, and we here focus on varieties of that form,~i.e.,  
\begin{equation}\label{eq:very_affine}
    X \,=\, \left(\CC^\ast\right)^n\setminus V\left( f_1\cdots f_\ell\right), \, 
\end{equation}    
where $(f_1,\ldots,f_\ell)$ is a family of Laurent polynomials.
In contrast to the setup of \eqref{eq:integral}, the coefficients of the $f_i$'s are now considered to be fixed. Via the graph morphism, 
$X$ as in~\eqref{eq:very_affine} embeds into the $(n+\ell)$-dimensional algebraic torus and hence is very affine.

For cohomological considerations, we denote by $\Omega^k_X$ the space of global algebraic differential $k$-forms on~$X$, which we assume to be smooth. Consider the complex $(\Omega_X^\bullet,\nabla_\omega)$, with the differential $\nabla_\omega^{(k)}\colon \Omega_X^k \to \Omega_X^{k+1}$ being twisted by the logarithmic 
one-form
\[ \omega \,\coloneqq\, \operatorname{dlog}(f^sx^\nu)\,=\, \sum_{j=1}^{\ell} s_j\frac{\d f_j}{f_j} + \sum_{i=1}^n \nu_i \frac{\d x_i}{x_i} \, ,\] i.e., $\nabla_\omega=\d \, + \omega \wedge \, .$ The cohomology vector space \[ H^k(X,\omega) \,\coloneqq \, \ker \big(\nabla_\omega ^{(k)}\big)/\operatorname{im}  \big(\nabla_\omega ^{(k-1)}\big)\] is the $k$-th {\em twisted cohomology space} of~$X$.
It is generated by equivalence classes of elements of the form $f^ax^b \, \frac{\d x_{i_1}\wedge \cdots \wedge \d x_{i_k}}{x_{i_1}\cdots \,  x_{i_k}}$ with integer vectors~$a,b$. 

Since $s,\nu$ can be vectors of complex numbers, one needs to consider a branch cut of the multi-valued function to make sense of the integral in~\eqref{eq:integral}. This is what is actually being taken care of by the twist by~$\omega$. The way this issue is taken into account is more explicit in the description of twisted homology. 
In fact, an element in the $k$-th twisted homology vector space of $X$ can be thought of as a 
singular chain on~$X$ of dimension~$k$ together with the choice of a branch cut of~$f^sx^\nu$. For elaborations on twisted homology, we refer to~\cite{AFST22}. In this presentation, we define the twisted homology spaces $H_k(X,\omega)$ simply as the dual vector spaces of the $H^k(X,\omega)$.
%
%
Despite the complexity in the construction, a drastic simplification happens when studying these twisted (co)-homology vector spaces, due to a vanishing theorem for all the intermediate cohomology spaces.

For fixed $(s,\nu)$, one has the perfect pairing 
\begin{align}\begin{split}\label{eq:pairing}
    \langle \cdot,\cdot \rangle \colon \, H^n(X,\omega) \times H_n(X,\omega) \, \longrightarrow \, \CC \, ,\\
    \left([\eta],[\Gamma]\right)\mapsto \int_\Gamma f^s x^{\nu} \, \eta  \qquad 
\end{split}\end{align}
between the $n$-th twisted cohomology and homo\-lo\-gy spaces of~$X$. When the parameters are rational numbers, numbers that can be expressed as an integral of the form~\eqref{eq:pairing} are {\em periods} of the variety~$X$. Recent progress in high precision evaluation of such periods to thousands of digits is discussed in~\cite{LPPV24}.

Let $s_j,\nu_i\in \CC$ be fixed generic complex numbers, then the $\CC$-vector space spanned by the functions
\begin{align}\label{eq:Iab} 
\Big\{ \left[ \Gamma \right] \,\mapsto\, \cI_{a,b} \coloneqq \int_\Gamma f^{s+a}x^{\nu+b} \, \frac{\d x}{x} \,  \Big\}_{(a,b)\in \ZZ^\ell \times \ZZ^n}
\end{align}
when allowing integer shifts of $s$ and $\nu$, is finite-dimensional. Its dimension equals the signed Euler characteristic of $X$, which coincides with the dimension of $H^n(X,\omega)$. 
An overview of these results is presented in~\cite{AFST22}. The use of twisted de Rham 
cohomology for Feynman integrals was initiated in~\cite{MastroliaMizera}. If the integral in~\eqref{eq:Iab} comes from a Feynman diagram, then the integrals in the basis are called ``master integrals.'' Describing a set of master integrals is extremely important in the study of scattering amplitudes. In this way, the evaluation of Feynman integrals reduces to the evaluation of master integrals, which significantly reduces the complexity.

\subsection{Algebraic analysis}\label{sec:algana}
Algebraic analysis is a mathematical field that investigates linear partial differential equations (PDEs) with methods from algebraic geometry, category theory, and complex analysis.
Homogeneous, linear PDEs with polynomial coefficients are encoded as elements of the ($n$-th) {\em Weyl algebra}, denoted 
$$ D_n \,=\, \CC[x_1,\ldots,x_n]\left\langle \partial_1,\ldots,\partial_n \right\rangle,$$
or just $D$,
which encodes linear differential operators with polynomial coefficients. Allowing coefficients in rational functions 
$\CC(x_1,\ldots,x_n)$,  results in the {\em rational Weyl algebra}, denoted $R_n$. The Weyl algebra is non-commutative: all of its generators are assumed to commute, except  $\partial_i$ and $x_i$, $i=1,\ldots,n$. They obey Leibniz' rule, i.e., their commutator $[\partial_i,x_i]=\partial_ix_i-x_i\partial_i$ fulfills $[\partial_i,x_i]=1$. 
We denote the action of a differential operator on a function $f(x_1,\ldots,x_n)$ by a bullet; e.g., $\partial_i\bullet f=\frac{\partial f}{\partial x_i}$. 

A system of PDEs corresponds to a left ideal $I\subset D$ in the Weyl algebra, to which one associates the \mbox{$D$-module}~$D/I$. Among them, ``holonomic'' $D$-ideals or \mbox{$D$-modules} give rise to finite-dimensional solution spaces. Also concepts like singularities (in the sense of poles or essential singularities) and integral transforms, such as the Fourier--Laplace or Mellin transform of functions, have their algebraic counterpart in this language. Even the growth behavior of solutions to the considered system of PDEs is captured by the corresponding $D$-ideal: functions that have moderate growth when approaching singularities give rise to {\em regular} singularities of the $D$-module, and to {\em irregular} singularities otherwise.

One may, however, not simply switch between a \mbox{$D$-module} and its solution space. For steps towards deeper structural insights of that kind, one needs to pass on to a more abstract setup, namely to sheaves of $\mathcal{D}_X$-modules on a smooth algebraic variety~$X$. In this generalized setup, one can phrase a positive answer to Hilbert's $21$st problem, which has a negative answer in its original formulation. 
Holomorphic solutions of a $\mathcal{D}_X$-module $\mathcal{M}$ are recovered algebraically as the $\mathcal{D}_X$-linear morphisms $\mathcal{H}\!\operatorname{om}_{\mathcal{D}_X}(\mathcal{M},\mathcal{O}_X^{\text{an}})$ from $\mathcal{M}$ to holomorphic functions $\mathcal{O}_X^{\text{an}}$ on~$X$. The right-derived functor of $\mathcal{H}\!\operatorname{om}_{\mathcal{D}_X}(\cdot,\mathcal{O}_X^{\text{an}})$ is called {\em solution functor} and implies a categorical equivalence of regular holonomic $\mathcal{D}_X$-modules and their topological counterpart: their solution complexes. This is a consequence of the Riemann--Hilbert correspondence, see e.g.~\cite{HTT08}.

In applications, one works mostly over the affine space. The $n$-th Weyl algebra $D$ is precisely the ring of global sections of $\mathcal{D}_{\CC^n}$, and $D$-modules are in one-to-one correspondence with those sheaves of $\mathcal{D}_{\CC^n}$-modules that are quasi-coherent over $\mathcal{O}_{\CC^n}$. 
In the holonomic case, the $D$-module structure then corresponds to an integrable connection on the underlying vector bundle. After trivializing the bundle, the action of the Weyl algebra can be written in connection form $\mathrm{d}+M\wedge \ $, with $M$ a matrix of differential one-forms. Dually, via the tensor-hom adjunction, the connection matrix can equivalently be expressed in the tangent space. The resulting system of PDEs hence gets expressed in matrix form, generalizing the companion matrix of ODEs, and can be computed 
by using Gröbner bases in the~Weyl~algebra~\cite{SST00}.
\begin{example}\label{ex:connection}
Let $I$ be the $D_2$-ideal generated by $P_1 =x_2^2\partial_2^2+(x_1+2x_2)\partial_{2}, $ $ P_2 = x_1{x_2}\partial_1\partial_{2}-(x_1+{x_2})\partial_{2}$, and $ 
P_3 =x_1^2\partial_{1}^2+x_1\partial_{1}+(x_1+{x_2})\partial_{2}.$
For the degree reverse lexicographic order on $R_2$,
the set of standard monomials of $R_2I$ is $S=\{ 1,\partial_2,\partial_1 \}$. 
For $f\in \operatorname{Sol}(I)$, denote $F=(f,\partial_2\bullet f,\partial_1\bullet f)^\top$. In the $\CC(x_1,x_2)$-basis~$S$ of ${R_2}/{R_2I}$, the connection matrices of $I$ are
\begin{align*}
    M_1 = \left[\begin{smallmatrix} 
   0 & 0 & 1\\
   0  & \frac{x+y}{xy} & 0\\
   0 & -\frac{x+y}{x^2} & -\frac{1}{x}\end{smallmatrix}\right] \text{ and } \    M_2 = \left[\begin{smallmatrix} 0 & 1 & 0 \\
   0 & -\frac{x+2y}{y^2} & 0\\
   0 & \frac{x+y}{xy} & 0
   \end{smallmatrix}\right],
\end{align*}
so that 
$\partial_i \bullet F = M_i \cdot F$, $i=1,2$, 
for any $f\in \operatorname{Sol}(I)$.
\end{example}

The {\em holonomic rank} of a $D_n$-ideal $I$  is defined as
$$ \operatorname{rank}(I) \,=\, \dim_{\CC(x_1,\ldots,x_n)} \left( R_n / R_n I \right) .$$
For instance, the $D_2$-ideal $I$ in \Cref{ex:connection} 
 has holonomic rank~$3$.
On simply connected domains in~$\CC^n$ outside the singular locus of a $D$-ideal $I$, the dimension of the space of holomorphic solutions is given by the holonomic rank of the $D$-ideal. This follows from the theorem of Cauchy, Kovalevskaya, and Kashiwara. 
As a consequence, functions $f$ that are solutions to holonomic $D$-ideals, called ``holonomic functions,'' can be encoded by finite data as follows. Given a subideal $I\subset \operatorname{Ann}_D(f)$ of finite holonomic rank, where $\operatorname{Ann}_D(f)$ denotes the annihilating \mbox{$D$-ideal} of~$f$,
$$ \operatorname{Ann}_D(f) \,=\, \{ P
\in D \,|\, P\bullet f =0\}\, ,$$
the function $f$ is uniquely determined by $I$ while prescribing $\operatorname{rank}(I)$-many initial conditions. 
Various functions throughout the sciences are holonomic. Such functions can be encoded, manipulated, and evaluated in terms of symbolic computations of their annihilating $D$-ideal. Holonomic functions were first treated algorithmically by Zeilberger. 

Ideals in $D$ encode crucial properties of their solution functions. For instance, $f(x_1,\ldots,x_n)$ being a solution~to $\theta_1 + \cdots + \theta_n-k \in D_n$,  
where $\theta_i=x_i\partial_i$ denotes the $i$th {Euler operator} and $k\in \NN$, implies that $f$ is homogeneous of degree $k$, 
i.e., $f(tx_1,\ldots,tx_n)=t^k\cdot f(x_1,\ldots,x_n)$ by the converse of Euler's homogeneous function theorem.

Prominent examples of $D$-ideals with a strong combinatorial flavor are GKZ systems~\cite{GKZ90}, $H_A(\kappa)$, also called ``$A$-hypergeometric systems.'' They are encoded by an integer matrix $A\in \ZZ^{k\times n}$, and a parameter vector $\kappa\in \CC^n$. The $D$-ideal $H_A(\kappa)$ is the $D$-ideal generated by the {\em toric ideal} of~$A$, 
$$I_A \,=\, \left\langle \, \partial^u-\partial^v \,|\, u-v\in\ker(A), \, u,v\in \NN^n\, \right\rangle,$$ 
and the entries of $A\cdot (\theta_1,\ldots,\theta_n)^\top-\kappa$.
Their solution functions are generalized Euler integrals. More precisely, in this context, the integral~\eqref{eq:integral} is viewed as a function of the coefficients $c$ of the Laurent polynomials~$f_i$. The GKZ system $H_A(\kappa)$ is therefore an ideal in the Weyl algebra in the $c$-variables. Feynman integrals in their Lee--Pomeransky representation are examples of generalized Euler integrals, where the coefficients of the graph polynomial are constrained to linear subspaces of the kinematic space---the latter being the affine space with coordinates given by the linearly independent kinematic variables. Therefore, they are solutions to restrictions of GKZ systems. However, computing restrictions of \mbox{$D$-ideals} is computationally challenging, which calls for the development of new strategies. Further challenges to tackle are non-generic parameter values for GKZ systems, which often constitute interesting and relevant cases. 

\smallskip

A concrete example of the utility of connecting Feynman integrals to GKZ systems lies in the Landau analysis, namely the study of the locus of complex parameters for which a Feynman integral develops singularities. In~\cite{FevolaMizeraTelen}, the singularities of a Feynman integral are formalized as the subspace of physical parameters on which the Feynman integrand is “more singular” than generically. More specifically, this is done by determining the locus of kinematics for which the signed Euler characteristic of the very affine variety fixed by the choice of coefficients is smaller than its generic value. This idea was inspired by the GKZ framework, in which the number of linearly independent solutions is bounded above by the signed Euler characteristic of the variety in~\eqref{eq:very_affine}, and the locus in the coefficient space for which this bound is not achieved characterizes the singularities of the generalized Euler integral. This locus is a hypersurface in the coefficient space and its defining polynomial is called ``principal $A$-determinant.'' However, when this polynomial is restricted to the subspace that is of physical interest, it often trivializes and fails to yield meaningful information about the singularities. To address this problem, \cite{FevolaMizeraTelen} proposes methods and algorithms to first specialize to the desired subspace and to then extract singularity information.

\subsection{Combinatorial algebraic geometry}\label{sec:combalggeo}
Various combinatorial aspects of algebraic geometry help to study scattering amplitudes and cosmological correlators. We here focus on positive geometries~\cite{PosGeom,Lam}, and start right away with giving their definition.

Let $X$ be a complex, irreducible, $d$-dimensional algebraic variety that is defined over the real numbers. Typically, $X$ is chosen to be projective. Equip its real-valued points, $X(\RR)$, with the analytic topology. Let $X_{\geq 0}$ be a closed, semi-algebraic subset of $X(\RR)$ such that its interior, $X_{>0}\coloneqq\operatorname{Int}(X_{\geq 0})$, is an oriented $d$-manifold, and such that the analytic closure of $X_{>0}$ in $X(\RR)$ recovers $X_{\geq 0}$. Denote by $\partial X_{\geq 0}=X_{\geq 0}\setminus X_{> 0} $ the boundary, by $\partial X$ the Zariski closure of $\partial X_{\geq 0}$ over~$\CC$, and $C_1,\ldots,C_r$ the irreducible components of $\partial X$. The {\em facets} of $X_{\geq 0}$ are the closures $C_{i,\geq 0}$ of the interiors of $C_i\cap X_{\geq 0}$ in~$C_i(\RR)$.

The tuple $(X,X_{\geq 0})$ is a {\em positive geometry} if there exists a unique, non-zero, rational\footnote{``Rational'' in this context means that the differential form is defined in terms of rational functions on the variety~$X$, i.e., elements of the function field of~$X$.}
differential $d$-form $\Omega(X,X_{\geq 0})$, called the {\em canonical form}, satisfying the following recursive axioms:
   If $d=0$, then 
    $X=X_{\geq 0}$ is a point, and one defines $\Omega(X,X_{\geq 0})=\pm 1$, depending on the orientation.
    If $d>0$, one requires that $\Omega(X,X_{\geq 0})$ has poles along the boundary components $C_i$ only, that these poles are at most simple, and that for each $i=1,\ldots,r$, the residue of $\Omega(X,X_{\geq 0})$ along $C_i$ is the canonical form of $(C,C_{i,\geq 0})$.

For instance, any closed interval $[a,b]\subset \RR$ defines a positive geometry $\big(\PP_{\CC}^1,[a,b]\big)$, whose canonical form is $\dlog\big(\frac{x-b}{x-a}\big)$ (up to a sign). Note that $(\PP_{\CC}^1,\PP_{\RR}^1)$ is not a positive geometry: $\partial \PP_{\RR}^1 = \emptyset $, but there is no global differential one-form on~$\PP_{\CC}^1$ without any poles.

\begin{example}
In $\PP_{\CC}^n$ with homogeneous coordinates $x_0,\ldots,x_n$, consider the projective simplex \mbox{$\Delta^n=\PP_{\geq 0}^n(\RR)$} consisting of points in $\PP_{\RR}^n$ that are representable by non-negative coordinates. The pair $\left(\PP_{\CC}^n,\Delta^n\right)$ is a positive geometry. In the affine chart \mbox{$\{x_0\neq 0\}$}, denoting $y_i=\frac{x_i}{x_0}$, its canonical form is
$$ \Omega \left( \PP_{\CC}^n,\Delta^n \right) \,=\, \wedge_{i=1}^n \frac{\d y_i}{y_i} \,=\, \wedge_{i=1}^n \dlog y_i \, . $$
\end{example}

More generally, also polytopes are positive geometries. Due to the additivity of the canonical form with respect to triangulations, it is sufficient to know the canonical forms of simplices in order to compute that of polytopes. Alternatively, it can obtained in terms of the adjoint of the polytope.

Normal, projective, toric varieties $X_P$, with $P$ a lattice polytope, constitute further instances of positive geometries: they have a natural, positive part $X_{P,>0}\subset X_P$ sitting inside them.

An example of a positive geometry, which is prominent in physics, is the positive Grassmannian, whose construction we recall now. Let $\mathbb{K}$ be a field. For natural numbers $k,n\in \NN$, the set of $k$-dimensional subspaces of $\mathbb{K}^n$, denoted $\Gr_{k,n}(\mathbb{K})$, is an algebraic variety in a natural way, called {\em Grassmannian}. Any point $V$ of the Grassmannian can be represented as the row-span of a $k\times n$ matrix $M_V$ of full rank, with entries in $\mathbb{K}$. Let $\binom{[n]}{k}$ denote the set of all $k$-element subsets of $[n]=\{1,\ldots,n\}$. For $I\in \binom{[n]}{k}$, the corresponding $k\times k$ minor of~$M_V$ is denoted by~$p_I(V)$. These minors are called {\em Plücker coordinates} of~$V$. Up to a simultaneous rescaling, they are independent of the choice of the representative~$M_V$, and they fulfill homogeneous, quadratic relations, called ``Plücker relations.'' Grassmannians hence are projective varieties via the {\em Plücker embedding}
\begin{align*}
\Gr_{k,n} \longrightarrow \mathbb{P}^{\binom{n}{k}-1} ,\  V\mapsto (p_I(V))_{I\in  \binom{[n]}{k}} 
\end{align*}
into projective space. 
The {\em non-negative Grassmannian} (informally called {{\em positive Grassmannian}), denoted $\Gr_{k,n}^{\geq 0}$, is the semi-algebraic subset of 
$\Gr_{k,n}(\mathbb{R})$, where all Plücker coordinates are non-negative.

\begin{example}
Consider the matrix
$M_V=\left[\begin{smallmatrix}
1 & 0 & -2 & -6\\
0 &1 &1 &3
\end{smallmatrix}\right]$ 
of full rank. Its row span corresponds to a point $V\in\Gr_{2,4}(\RR)$. Since $p_{12}(M_V) = 1,$ $ p_{13}(M_V) = 1,$ $ p_{14}(M_V) = 3, $ $ p_{23}(M_V) = 2,$ $ p_{24}(M_V) = 6,$ and $p_{34}(M_V) = 0$, $M_V$ represents a point in~$\Gr_{2,4}^{\geq 0}$.
\end{example}

By \cite[Theorem 9]{Lam}, the pair $(\Gr_{k,n}, \Gr_{k,n}^{\geq 0})$ is a simplex-like positive geometry. Its faces are called {\em positroid cells}, and their Zariski closures are {\em positroid varieties}. Images of positive Grassmannians under ``positive'' linear maps are called {\em amplituhedra}, as introduced by Arkani-Hamed and Trnka~\cite{Amplitu}. To be precise, the {amplituhedron} $\mathcal{A}_{k,n,m}(Z)\subset \Gr_{k,k+m}$ is the image of $\Gr_{k,n}^{\geq 0}$ under a linear map that is induced by right multiplication by a real $n\times (k+m)$ matrix~$Z$, all of whose maximal minors are positive. It is conjectured that also  amplituhedra are positive~geometries.

In their recent work~\cite{BrownDupont}, Brown and Dupont define logarithmic canonical forms using relative homology and mixed Hodge theory.

\section{Interactions}\label{sec:interactions}
\vspace*{-2mm}
In this last section, we showcase a small assortment of the main problems and questions driving some of the recent research directions in this interplay between algebra, combinatorics, geometry, and physics.

In order to deal with divergent integrals, dimensional regularization allows the spacetime dimension $D=4-2\varepsilon$ to be non-integer. In this framework, a common strategy in particle physics is to  construct matrix differential equations behind Feynman integrals in ``$\varepsilon$-factorized form,'' i.e., of the form {$\d 
{\cI}=\varepsilon M\cI $,} with ${\cI}$ a vector of master integrals, see~Section~\ref{sec:alggeo}.
This special form is beneficial, since for such systems, Laurent series solutions in~$\varepsilon$ can be computed via the path-ordered exponential formalism.
In~\cite{HPSZ23}, Gröbner basis techniques in the Weyl algebra are used to simplify and solve linear PDEs behind Feynman integrals. 
PDEs in cosmology are tackled in \cite{DEcosmological,CosmoDmod}, where the $\varepsilon$-factorized form of the differential equation was made explicit. The structure of the differential equations naturally carries a strong combinatorial flavor, leading to a variety of intriguing open questions. 

The drive to improve techniques for deriving differential equations stems from the computational challenges posed by physics. As the complexity of fundamental interactions increases—such as when considering graphs with more cycles and intricate topologies—the number of parameters grows significantly. Additionally, the exceptional nature of the parameters in physics suggests to focus also on non-generic cases. This effort employs techniques from nonlinear algebra \cite{NonlinearAlgebra}, which uses both symbolic and numerical methods, like Gröbner basis computations and homotopy continuation. Various software and programming languages, including {\tt Julia}, {\tt Macaulay2}, {\tt Maple}, {\tt Mathematica}, {\tt Polymake}, {\tt Sage}, and {\tt Singular}, offer effective tools for these computations.

The interplay between symbolic and numerical techniques is a key strength linking the problem of maximum likelihood estimation in algebraic statistics to the study of scattering equations in particle physics, as first established in~\cite{ST21}. Both fields involve maximizing a logarithmic function and counting the solutions to a system of rational equations over a very affine variety. For instance, the likelihood equations for the moduli space $\mathcal{M}_{0,n}$ of $n$ marked points on a curve of genus $0$---when interpreted as a statistical model---are the scattering equations for an $n$-point tree-level scattering process. This connection allows methods from likelihood geometry to be applied to the investigation of scattering amplitudes. Moreover, the development of degeneration techniques for solving such systems is an emerging area of research.

Positive Grassmannians 
serve as the space of kinematic parameters. Their images under positive linear maps, namely amplituhedra, encode properties of scattering amplitudes in planar $\mathcal{N}=4$ SYM 
theory and drastically simplify the computation of gluon amplitudes. The scattering amplitude is conjecturally obtained by the canonical form of the amplituhedron at tree-level, and as the integral over a regulator function against the canonical form of the amplituhedron at loop-level, see \cite[p.~17]{Amplitu} for the definition of the loop amplituhedron. In practice, the canonical form is computed via subdivisions into positroid cells, referred to as ``BCFW triangulations.'' The goal of proving the positive nature of the amplituhedron opens a research area towards understanding, 
e.g., subdivisions and boundaries of the amplituhedron with ideas being inspired from the theory of matroids, tropical geometry, and cluster algebras; see~\cite{Williams} and the references therein. Current efforts aim to extend the use of positive geometries to further~QFTs.

We intended to provide a glimpse into the main research areas within the context we are discussing, but the connections and questions are wide-ranging. Among these, we stress the link between moduli spaces and quantum field theory, non-vanishing classes of cohomology of graph complexes and Feynman integrals, cluster and flag varieties in the context of spinor-helicity varieties, as well as connection and Macaulay matrices.
In short, physics suggests new, intriguing mathematical structures. It is now left to the community to catch up with nailing down the details of the emerging mathematical objects and theories, and to attest them. As verified by various fruitful and successful collaborations, important first steps have already been taken.

\bigskip 
\noindent
{\bf Authors' addresses:}

\smallskip

\noindent Claudia Fevola,
Université Paris-Saclay, Inria, Palaiseau, France 
\\{\tt claudia.fevola@inria.fr}

\medskip
\noindent Anna-Laura Sattelberger, Max Planck Institute for Mathematics in the Sciences, Leipzig, Germany
\\{\tt anna-laura.sattelberger@mis.mpg.de}

\vspace*{-1mm}
\paragraph{Funding statement}
CF has received funding from the
European Union’s Horizon 2020 research and innovation programme under the Marie Sk\l odowska-Curie grant agreement No 101034255.
The research of ALS is funded by the European Union (ERC, UNIVERSE PLUS, 101118787). Views and opinions expressed are, however, those of the author(s) only and do not necessarily reflect those of the European Union or the European Research Council Executive Agency. Neither the European Union nor the granting authority can be held responsible for them.

\end{document}